\documentclass[11pt]{amsart}
\usepackage{amsmath,amssymb}
\newtheorem{theorem}{Theorem}[section]
\newtheorem{corollary}[theorem]{Corollary}

\newtheorem{proposition}[theorem]{Proposition}
\newtheorem*{Y}{Yamabe Conjecture}
\newtheorem*{C1}{Compactness Conjecture}
\newtheorem*{C2}{Conjecture}
\newtheorem*{C3}{Weyl Vanishing Conjecture}

\begin{document}

\title[Yamabe problem]{Recent progress on the Yamabe problem}
\author{Simon Brendle and Fernando C. Marques}
\address{Department of Mathematics \\ Stanford University \\ 450 Serra Mall, Bldg 380 \\ Stanford, CA 94305} 
\address{Instituto de Matem\'atica Pura e Aplicada (IMPA) \\ Estrada Dona Castorina 110 \\ 22460-320 Rio de Janeiro \\ Brazil}
\thanks{The first author was supported in part by the National Science Foundation under grant DMS-0905628. The second author was supported by CNPq-Brazil, FAPERJ, and the Stanford Department of Mathematics.}
\dedicatory{Dedicated to Professor Richard Schoen on the occasion of his sixtieth birthday.}
\begin{abstract}
We give a survey of various compactness and non-compactness results for the Yamabe equation. We also discuss a conjecture of Hamilton concerning the asymptotic behavior of the parabolic Yamabe flow.
\end{abstract}
\maketitle

\section{The Yamabe problem}

A basic question in differential geometry is to find canonical metrics on a given manifold $M$. For example, if $\dim M = 2$, the uniformization theorem guarantees the existence of a metric of constant Gaussian curvature in any given conformal class:

\begin{theorem}
Let $M$ be a compact surface, and let $g_0$ be a Riemannian metric on $M$. Then there exists a smooth function $w: M \to \mathbb{R}$ with the property that the metric $g = e^{2w} \, g_0$ has constant Gaussian curvature.
\end{theorem}

Recall that two metrics $g_0$ and $g$ are said to be conformal if $g = e^{2w} \, g_0$ for some smooth function $w: M \to \mathbb{R}$.

We next consider the case $\dim M \geq 3$. Motivated by the uniformization theorem, Yamabe \cite{Yamabe} proposed the following conjecture: 

\begin{Y}
Let $(M,g_0)$ be a compact Riemannian manifold of dimension $n \geq 3$. Then there exists a metric $g$ on $M$ which is conformal to $g_0$ and has constant scalar curvature.
\end{Y}

The Yamabe problem can be reduced to the solvability of a certain semilinear elliptic equation. To that end, let us write $g = u^{\frac{4}{n-2}} \, g_0$ for some positive function $u$. Then the scalar curvature of $g$ is given by 
\begin{equation} 
\label{scalar}
R_g = u^{-\frac{n+2}{n-2}} \, \Big ( -\frac{4(n-1)}{n-2} \, \Delta_{g_0} u + R_{g_0} \, u \Big ), 
\end{equation}
where $\Delta_{g_0}$ and $R_{g_0}$ denote the Laplace operator and scalar curvature of the metric $g_0$. Therefore, the metric $g$ has constant scalar curvature $c$ if and only if $u$ is a solution of the equation 
\begin{equation} 
\label{yamabe.equation}
\frac{4(n-1)}{n-2} \, \Delta_{g_0} u - R_{g_0} \, u + c \, u^{\frac{n+2}{n-2}} = 0. 
\end{equation}
The equation (\ref{yamabe.equation}) is referred to as the Yamabe equation. 

The Yamabe equation has a natural variational characterization. To describe this, we consider the normalized Einstein-Hilbert action 
\[\mathcal{E}(g) = \frac{\int_M R_g \, d\text{\rm vol}_g}{\text{\rm vol}(M,g)^{\frac{n-2}{n}}}.\] 
A metric $g$ is a critical point of $\mathcal{E}$ if and only if $g$ is an Einstein metric. We now restrict the functional $\mathcal{E}$ to the conformal class of $g_0$. Given any positive function $u$, the Yamabe functional is defined by 
\[E_{g_0}(u) = \mathcal{E}(u^{\frac{4}{n-2}} \, g_0).\] 
The identity (\ref{scalar}) implies 
\[E_{g_0}(u) = \frac{\int_M \big ( \frac{4(n-1)}{n-2} \, |du|_{g_0}^2 + R_{g_0} \, u^2 \big ) \, d\text{\rm vol}_{g_0}}{\big ( \int_M u^{\frac{2n}{n-2}} \, d\text{\rm vol}_{g_0} \big )^{\frac{n-2}{n}}}.\] 
A positive function $u$ is a critical point of the functional $E_{g_0}$ if and only if $u$ satisfies (\ref{yamabe.equation}) for some constant $c$. Therefore, the Yamabe problem is equivalent to finding critical points of the functional $E_{g_0}$.

The Yamabe constant of a compact manifold $(M,g_0)$ is defined as the infimum of the Yamabe functional $E_{g_0}$, i.e. 
\[Y(M,g_0) = \inf_{0 < u \in C^\infty(M)} E_{g_0}(u).\] 
To fix notation, we denote by $Y(S^n)$ the Yamabe constant of the sphere $S^n$, equipped with its standard metric.

The solution of the Yamabe problem involves two steps. First, Trudinger and Aubin gave a sufficient condition for the existence of a minimizer.
Notice that it is not difficult to show that $Y(M,g_0) \leq Y(S^n)$.

\begin{theorem}[T.~Aubin \cite{Aubin}; N.~Trudinger \cite{Trudinger}] 
\label{existence.of.minimizer}
Let $(M,g_0)$ be a compact Riemannian manifold with $Y(M,g_0) < Y(S^n)$. Then the infimum of the functional $E_{g_0}$ is attained. In particular, there exists a positive function $u$ and constant $c$ such that 
(\ref{yamabe.equation}) holds.
\end{theorem}

Second, it was shown by Aubin and Schoen that the inequality $Y(M,g_0) < Y(S^n)$ holds unless $(M,g_0)$ is conformally equivalent to the standard sphere.

\begin{theorem}[T.~Aubin \cite{Aubin}; R.~Schoen \cite{Schoen1}] 
\label{inequality.for.yamabe.constant}
Let $(M,g_0)$ be a compact Riemannian manifold of dimension $n \geq 3$ which is not conformally equivalent to the standard sphere $S^n$. Then $Y(M,g_0) < Y(S^n)$.
\end{theorem}

The proof of Theorem \ref{inequality.for.yamabe.constant} naturally breaks into two cases.  If $n \geq 6$ and $(M,g_0)$ is not conformally flat, Aubin \cite{Aubin} was able to construct a test function $u: M \to \mathbb{R}$ with $E_{g_0}(u) < Y(S^n)$. This construction is purely local, and exploits the Weyl tensor of $(M,g_0)$. On the other hand, if $3 \leq n \leq 5$ or $(M,g_0)$ is conformally flat, Schoen \cite{Schoen1} showed that there exists a function $u$ with $E_{g_0}(u) < Y(S^n)$. Schoen's argument is global in nature and uses the Positive Mass Theorem (cf. \cite{Schoen-Yau}, \cite{Witten}) in a crucial way. For more details, we refer the reader to the survey of Lee and Parker \cite{Lee-Parker}.

Theorems \ref{existence.of.minimizer} and \ref{inequality.for.yamabe.constant} imply that the Yamabe problem has a solution for every compact manifold $(M,g_0)$. Bahri \cite{Bahri} gave an alternative proof in the locally conformally flat case. The method in \cite{Bahri} produces a critical point of the Yamabe functional which may not be a minimizer.

If $Y(M,g_0) \leq 0$, the solution to the Yamabe problem is unique, up to scaling. This is no longer true if $Y(M,g_0)>0$.  For example, let us consider the cylinder $S^{n-1} \times S^1$. 

\textbf{Example.} Let $(M,g_0) = S^{n-1}(1) \times S^1(\ell)$, where $S^{n-1}(1)$ denotes the $(n-1)$-sphere equipped with its standard metric and the $S^1$ factor has length $\ell$. If $\ell > 0$ is sufficiently small, the product metric is the unique solution of the Yamabe equation on $(M,g_0)$, up to scaling. However, if we choose $\ell$ sufficiently large, the equation (\ref{yamabe.equation}) will have an arbitrarily large number of solutions. See \cite{Schoen2} for more details.

More generally, Pollack \cite{Pollack} has obtained the following non-uniqueness theorem:

\begin{theorem}[D.~Pollack \cite{Pollack}]
\label{nonuniqueness}
Let $(M,g_0)$ be a compact manifold with $Y(M,g_0) > 0$. Given any positive integer $N$, there exists a metric $g$ such that $\|g - g_0\|_{C^0} \leq \frac{1}{N}$ and the equation 
\[\Delta_g u - \frac{n-2}{4(n-1)} \, R_g \, u + n(n-2) \, u^{\frac{n+2}{n-2}} = 0\] 
has at least $N$ solutions.
\end{theorem}

\section{The Compactness Conjecture}

In light of Theorem \ref{nonuniqueness}, it is of interest to understand the space of all solutions to the Yamabe equation. In particular, one would like to develop a Morse theory for the Yamabe functional $E_{g_0}$. This is a delicate question, as the Yamabe functional does not satisfy the Palais-Smale condition. In a topics course at Stanford in 1988, Schoen conjectured that the set of solutions of the Yamabe equation is compact.

\begin{C1}[R.~Schoen \cite{Schoen3}, \cite{Schoen5}]
Let $(M,g_0)$ be a compact Riemannian manifold of dimension $n \geq 3$ which is not conformally equivalent to the standard sphere $S^n$. Moreover, let us fix a real number $c$. Then the set of all positive solutions of (\ref{yamabe.equation}) is compact in the $C^2$-topology.
\end{C1}

Note that the Compactness Conjecture is a key step towards developing a Morse theory for the Yamabe functional. This is discussed in more detail in \cite{Schoen3}.

The standard sphere $(S^n,\overline{g})$ is special because its group of conformal transformations is noncompact. By a theorem of Obata (see \cite{Obata}), a conformal metric $g = u^{\frac{4}{n-2}} \, \overline{g}$ has constant scalar curvature if and only if there exists a conformal transformation $\psi: (S^n,\overline{g}) \to (S^n,\overline{g})$ such that $g = \lambda \, \psi^*(\overline{g})$ for some positive constant $\lambda$. 

In order to approach the conjecture, Schoen proposed a strategy based on the so-called Pohozaev identity:

\begin{proposition}[Pohozaev Identity]
Let $(\Omega^n,g)$ be a Riemannian domain, $n \geq 3$. If $X$ is a vector field on $\Omega$, then
\[\frac{n-2}{2n} \int_\Omega X(R_g)\, d\text{\rm vol}_g + \int_\Omega \langle \mathcal{D}_g X,T_g \rangle \, d\text{\rm vol}_g = \int_{\partial \Omega} T_g(X,\eta_g)\, d\sigma_g.\] 
Here, $T_g = \text{\rm Ric}_g - \frac{1}{n} \, R_g \, g$ denotes the traceless Ricci tensor, $(\mathcal{D}_gX)_{ij}=X_{i;j}+X_{j;i} - \frac{2}{n} \, \text{\rm div}_g \, X \, g_{ij}$ is the conformal Killing operator, and $\eta_g$ is the outward unit normal to $\partial \Omega$.
\end{proposition}

The idea is to argue by contradiction. If $u_\nu$ is a nonconverging sequence of solutions to 
\begin{equation} 
\label{yamabe.equation.2}
\Delta_{g_0} u - \frac{n-2}{4(n-1)} \, R_{g_0} \, u + n(n-2) \, u^{\frac{n+2}{n-2}} = 0,
\end{equation}
then it follows from standard blow-up arguments that $u_\nu$ has to concentrate and form spherical bubbles at some points of the manifold. More precisely, there exist an integer $N_\nu > 0$ and a  finite set $\{y_{1,\nu}, \dots,y_{N_\nu,\nu}\} \subset M$ of local maxima of $u_\nu$ with the following properties: 
\begin{itemize}
\item For each $i=1,\dots,N_\nu$, the function $u_\nu$ is well approximated by the function
\[\Big ( \frac{\varepsilon_{i,\nu}}{\varepsilon_{i,\nu}^2 + d_{g_0}(y_{i,\nu},\cdot)^2} \Big )^{\frac{n-2}{2}}\] 
in a small neighborhood of $y_{i,\nu}$. Here, $\varepsilon_{i,\nu} = u_\nu(y_{i,\nu})^{-\frac{2}{n-2}}$. 
\item The function $u_\nu$ is bounded away from the set $\{y_{1,\nu}, \dots,y_{N_\nu,\nu}\}$. 
\end{itemize}
Notice that we are not assuming that the sequence $u_\nu$ has bounded energy.

First one has to deal with the case of isolated simple blow-up, in which there is no accumulation of more than one bubble at a single point. The general case of multiple blow-up can be reduced to the simple blow-up case using scaling arguments.

Suppose $\overline{y} = \lim_{\nu \rightarrow \infty} y_\nu$ is an isolated simple blow-up point. By assumption, the conformal metric $g_\nu = u_\nu^\frac{4}{n-2}g_0$ has constant scalar curvature. Applying the Pohozaev identity to the geodesic ball $B_\delta(y_\nu)=\{ p \in M: r = d_{g_0}(y_\nu,p) \leq \delta\}$, endowed with the Riemannian metric $g_\nu$ and the radial vector field $X = r \, \frac{\partial}{\partial r}$, $r = d_{g_0}(y_\nu,\cdot)$, yields
\begin{equation} 
\label{pohozaev_i}
\int_{B_\delta(y_\nu)} \langle \mathcal{D}_{g_\nu}X,T_{g_\nu} \rangle \, d\text{\rm vol}_{g_\nu} = \int_{\partial B_\delta(y_\nu)} T_{g_\nu}(X,\eta_{g_\nu})\, d\sigma_{g_\nu}.
\end{equation}

Schoen realized that the asymptotic behavior of the boundary integrals in (\ref{pohozaev_i}) is related to the expansion of the Green's function $G_L$ of the conformal Laplacian $L_{g_0} = \Delta_{g_0} - \frac{n-2}{4(n-1)} \, R_{g_0}$, with pole at $\overline{y}$. Notice that the manifold $M \setminus \{p\}$ equipped with the metric $\hat{g} = G_L^\frac{4}{n-2} \, g_0$ is asymptotically flat and has zero scalar curvature. If $n \leq 5$, the ADM mass $m$ of $\hat{g}$ is well defined. Since $(M,g_0)$ is not conformally equivalent to the standard sphere, the metric $\hat{g}$ is not flat, and hence $m > 0$ by the Positive Mass Theorem. After a calculation, it follows that
\[\lim_{\delta \rightarrow 0} \lim_{\nu \rightarrow \infty} \varepsilon_\nu^{2-n} \int_{\partial B_\delta(y_\nu)} T_{g_\nu}(X,\eta_{g_\nu})\, d\sigma_{g_\nu} < 0.\]
In order to get a contradiction, one would like to show that the corresponding limit of the left hand side of (\ref{pohozaev_i}) is nonnegative. When $(M,g_0)$ is locally conformally flat, it is possible to choose coordinates so that $r \, \frac{\partial}{\partial r}$ is a conformal Killing vector field, and therefore the left-hand side of (\ref{pohozaev_i}) is zero. A similar argument works in the three-dimensional case: in this case, basic estimates from above for $u_\nu$ imply that the left-hand side of (\ref{pohozaev_i}) converges to zero. These were the cases of the Compactness Conjecture covered by Schoen in the Stanford notes (see also \cite{Li-Zhu} and \cite{Schoen-Zhang}). In dimensions $4$ and $5$, the conjecture was proved by O.~Druet (see \cite{Druet}). 

The case $n \geq 6$ is more subtle. In this case, the manifold $(M \setminus \{p\},\hat{g})$ may not have a well-defined ADM mass unless the Weyl tensor of $g_0$ vanishes to order $[\frac{n-6}{2}]$ at the point $\overline{y}$. This causes extra difficulty, and it motivated Schoen to propose a conjecture concerning the location of possible blow-up points:

\begin{C3} [R.~Schoen \cite{Schoen5}] If $\overline{y} \in M$ is a blow-up point of a sequence of solutions $g_\nu = u_\nu^\frac{4}{n-2} \, g_0$ to the Yamabe Problem, then the Weyl tensor $W_{g_0}$ of the metric $g_0$ vanishes to order $[\frac{n-6}{2}]$ at the point $\overline{y}$. In other words, 
\[\limsup_{y \to \overline{y}} d_{g_0}(\overline{y},y)^{2-d} \, |W_{g_0}(y)| = 0,\] 
where $d = [\frac{n-2}{2}]$.
\end{C3}

The Weyl Vanishing Conjecture implies that the asymptotically flat manifold $(M \setminus \{p\},\hat{g})$ has a well-defined ADM mass. Furthermore, if $n \leq 7$, then the Positive Mass Theorem of Schoen and Yau \cite{Schoen-Yau} guarantees that the mass is positive. In \cite{Witten}, E. Witten proved the Positive Mass Theorem for spin manifolds, regardless of the dimension. Note that the Positive Mass Theorem in the locally conformally flat case was handled by a special argument in \cite{Schoen-Yau2}.

The Compactness and Weyl Vanishing Conjectures have been studied by numerous authors. It follows from the works of the second author and of Y.Y.~Li and L.~Zhang that compactness holds for all manifolds of dimension $n \leq 7$ (cf. \cite{Marques1}, \cite{Li-Zhang1}). Furthermore, Li and Zhang were able to verify the Compactness Conjecture for spin manifolds of dimension $n \leq 11$. The spin assumption is only needed in order to apply the Positive Mass Theorem.

The proof of the following compactness result does not require the Positive Mass Theorem, since there cannot be any blow-up points for purely local reasons. 

\begin{theorem}[Y.Y. Li and L. Zhang \cite{Li-Zhang1}, F.C. Marques \cite{Marques1}]
Let $(M,g_0)$ be a compact Riemannian manifold of dimension $n \geq 8$, not conformally diffeomorphic to the standard sphere. Suppose that 
\[|W_{g_0}(y)| + |\nabla W_{g_0}(y)| > 0\]
for all points $y \in M$. Then the set of solutions of (\ref{yamabe.equation.2}) is compact in the $C^2$-topology.
\end{theorem}

In June 2006, the first author constructed counterexamples to the Compactness Conjecture in dimension $n \geq 52$ (cf. \cite{Brendle4}). In \cite{Brendle-Marques}, the authors generalized this construction to dimension $25 \leq n \leq 51$. This construction will be discussed in more detail in the following section. Later, M. Khuri, the second author, and R. Schoen established the Compactness Conjecture for spin manifolds of dimension $n \leq 24$. In particular, the restriction on the dimension is sharp. Furthermore, the spin assumption can be removed if the Positive Mass Theorem holds true in general. We will give a detailed discussion of this compactness result in Section \ref{compactness}.

\section{Non-compactness results in dimension $n \geq 25$} 

\label{noncompactness}

In this section, we discuss counterexamples to Schoen's Compactness Conjecture. The following result was proved by the first author in \cite{Brendle4}:

\begin{theorem}[S.~Brendle \cite{Brendle4}]
\label{blow.up}
Assume that $n \geq 52$. Then there exists a smooth Riemannian metric $g$ on $S^n$ and a sequence of positive functions $\{v_\nu: \nu \in \mathbb{N}\}$ with the following properties: 
\begin{itemize}
\item $g$ is not conformally flat. 
\item For each $\nu \in \mathbb{N}$, the function $v_\nu$ is a solution of the Yamabe equation 
\[\Delta_g v_\nu - \frac{n-2}{4(n-1)} \, R_g \, v_\nu + n(n-2) \, v_\nu^{\frac{n+2}{n-2}} = 0.\] 
\item For each $\nu \in \mathbb{N}$, we have $E_g(v_\nu) < Y(S^n)$. Moreover, $E_g(v_\nu) \to Y(S^n)$ as $\nu \to \infty$. 
\item $\sup_{S^n} v_\nu \to \infty$ as $\nu \to \infty$.
\end{itemize}
\end{theorem}

We note that A.~Ambrosetti and A.~Malchiodi \cite{Berti-Malchiodi} have constructed examples of non-smooth background metrics on $S^n$ for which the set of solutions to the Yamabe equation is non-compact (see also \cite{Berti-Malchiodi} for related work). Moreover, it was shown by O.~Druet and E.~Hebey \cite{Druet-Hebey1}, \cite{Druet-Hebey2} that blow-up can occur for equations of the form $Lu + c \, u^{\frac{n+2}{n-2}} = 0$, where $L$ is a lower order perturbation of the conformal Laplacian on $S^n$.

We now describe the main ideas involved in the proof of Theorem \ref{blow.up}. We will consider small perturbations of the standard metric on $S^n$. It will be convenient to identify $\mathbb{R}^n$ with the complement of the north pole in $S^n$ via stereographic projection. This allows us to work on $\mathbb{R}^n$ instead of $S^n$.

Let us fix an integer $n \geq 52$. Moreover, we fix a constant tensor $W_{ijkl}$ which has the symmetries of the Weyl tensor. In particular, 
\[W_{ijkl} = -W_{jikl} = W_{klij}\] 
and 
\[\sum_{l=1}^n W_{ilkl} = 0.\] 
We assume that $W$ is non-trivial in the sense that $\sum_{i,j,k,l=1}^n (W_{ijkl} + W_{ilkj})^2 > 0$.

The following result is the key technical ingredient in the proof of Theorem \ref{blow.up}.

\begin{theorem}[\cite{Brendle4}, Proposition 24]
\label{gluing}
Assume that $\lambda$, $\mu$, and $\rho$ are real numbers satisfying $\mu \leq 1$ and $\lambda \leq \rho \leq 1$. Moreover, suppose that $h_{ik}(x)$ is a trace-free symmetric two-tensor on $\mathbb{R}^n$ with the following properties: 
\begin{itemize}
\item $h_{ik}(x) = \mu \, (\lambda^2 - |x|^2) \, \sum_{j,l=1}^n W_{ijkl} \, x_j \, x_l$ for all points $x \in \mathbb{R}^n$ with $|x| \leq \rho$.
\item $h_{ik}(x) = 0$ for all points $x \in \mathbb{R}^n$ with $|x| \geq 1$.
\item $\|h\|_{C^2} \leq \alpha$.
\end{itemize} 
Let us define a Riemannian metric $g_{ik}(x)$ on $\mathbb{R}^n$ by $g(x) = \exp(h(x))$. Then, if $\alpha$ and $\rho^{2-n} \, \mu^{-2} \, \lambda^{n-10}$ are sufficiently small, there exists a nonnegative function $v$ on $\mathbb{R}^n$ satisfying 
\[\Delta_g v - \frac{n-2}{4(n-1)} \, R_g \, v + n(n-2) \, v^{\frac{n+2}{n-2}} = 0\] 
and 
\[\sup_{|x| \leq \lambda} v(x) \geq c(n) \, \lambda^{\frac{2-n}{2}}.\] 
Here, $c(n)$ is positive constant that depends only on $n$.
\end{theorem}

In particular, if we choose $\lambda$ small enough, then the Yamabe equation will have a solution which is concentrated near the origin.

The proof of Theorem \ref{gluing} uses perturbation techniques. We will sketch the main ideas. Given any pair $(\xi,\varepsilon) \in \mathbb{R}^n \times (0,\infty)$, we consider the function 
\[u_{(\xi,\varepsilon)}(x) = \Big ( \frac{\varepsilon}{\varepsilon^2 + |x - \xi|^2} \Big )^{\frac{n-2}{2}}.\] 
The function $u_{(\xi,\varepsilon)}$ satisfies the differential equation 
\[\Delta u_{(\xi,\varepsilon)} + n(n-2) \, u_{(\xi,\varepsilon)}^{\frac{n+2}{n-2}} = 0.\] 
Using the implicit function theorem, one can construct a function $v_{(\xi,\varepsilon)}$ which is close to $u_{(\xi,\varepsilon)}$ in the $H^1$-topology and satisfies the differential equation 
\begin{align*} 
&\Delta_g v_{(\xi,\varepsilon)} - \frac{n-2}{4(n-1)} \, R_g \, v_{(\xi,\varepsilon)} + n(n-2) \, |v_{(\xi,\varepsilon)}|^{\frac{4}{n-2}} \, v_{(\xi,\varepsilon)} \\ 
&= \Big ( \frac{\varepsilon}{\varepsilon^2 + |x - \xi|^2} \Big )^{\frac{n+2}{2}} \, \bigg ( a_0(\xi,\varepsilon) \, \frac{\varepsilon^2 - |x - \xi|^2}{\varepsilon^2 + |x - \xi|^2 } + \sum_{i=1}^n a_i(\xi,\varepsilon) \, \frac{2\varepsilon \, (x_i - \xi_i)}{\varepsilon^2 + |x - \xi|^2} \bigg ).
\end{align*} 
for suitable coefficients $a_i(\xi,\varepsilon)$. For abbreviation, let 
\begin{align*} 
\mathcal{F}_g(\xi,\varepsilon) 
&= \int_{\mathbb{R}^n} \Big ( |dv_{(\xi,\varepsilon)}|^2 + \frac{n-2}{4(n-1)} \, R_g \, v_{(\xi,\varepsilon)}^2 - (n-2)^2 \, |v_{(\xi,\varepsilon)}|^{\frac{2n}{n-2}} \Big ) \\ 
&- 2(n-2) \, \Big ( \frac{Y(S^n)}{4n(n-1)} \Big )^{\frac{n}{2}}. 
\end{align*} 
We then have the following result: 

\begin{proposition}[\cite{Brendle4}, Proposition 6]
\label{lyapunov.schmidt}
For each point $(\bar{\xi},\bar{\varepsilon}) \in \mathbb{R}^n \times (0,\infty)$ the following statements are equivalent: 
\begin{itemize} 
\item The function $v_{(\xi,\varepsilon)}$ satisfies the differential equation 
\[\Delta_g v_{(\xi,\varepsilon)} - \frac{n-2}{4(n-1)} \, R_g \, v_{(\xi,\varepsilon)} + n(n-2) \, |v_{(\xi,\varepsilon)}|^{\frac{4}{n-2}} \, v_{(\xi,\varepsilon)} = 0.\] 
\item $a_i(\bar{\xi},\bar{\varepsilon}) = 0$ for $i = 0,1,\hdots,n$.
\item The point $(\bar{\xi},\bar{\varepsilon})$ is a critical point of the function $(\xi,\varepsilon) \mapsto \mathcal{F}_g(\xi,\varepsilon)$.
\end{itemize}
\end{proposition}

Therefore, the problem is reduced to finding critical points of the function $\mathcal{F}_g(\xi,\varepsilon)$. Unfortunately, we do not have an explicit expression for the function $\mathcal{F}_g(\xi,\varepsilon)$, as $v_{(\xi,\varepsilon)}$ is constructed using the implicit function theorem. To overcome this problem, we approximate the function $\mathcal{F}_g(\xi,\varepsilon)$ by another function which is explicitely computable. The following result is a direct consequence of Corollary 14 in \cite{Brendle4}.

\begin{proposition}[\cite{Brendle4}] 
\label{approximation.for.script.F}
Let $h_{ik}(x)$ be a trace-free symmetric two-tensor which satisfies the assumptions of Theorem \ref{gluing}. Moreover, let $g_{ik}(x)$ be the metric on $\mathbb{R}^n$ defined by $g(x) = \exp(h(x))$. Then 
\begin{align*} 
&|\lambda^{-8} \, \mu^{-2} \, \mathcal{F}_g(\lambda \xi,\lambda \varepsilon) - F(\xi,\varepsilon)| \\ 
&\leq C \, \lambda^{\frac{16}{n-2}} \, \mu^{\frac{4}{n-2}} + C \, \rho^{\frac{2-n}{2}} \, \mu^{-1} \, \lambda^{n-10} + C \, \rho^{2-n} \, \mu^{-2} \, \lambda^{n-10}
\end{align*} 
for all points $(\xi,\varepsilon)$ satisfying $|\xi| < 1$ and $\frac{n-8}{3(n+4)} < \varepsilon^2 < \frac{2(n-8)}{3(n+4)}$. Here, $F(\xi,\varepsilon)$ is a specific function which depends only on our choice of $W_{ijkl}$.
\end{proposition} 

It turns out that 
\[F(0,\varepsilon) = -I(n) \, \sum_{i,j,k,l=1}^n (W_{ijkl} + W_{ilkj})^2 \, \Big ( \frac{n-8}{n+4} \, \varepsilon^4 - 2 \, \varepsilon^6 + \frac{n+8}{n-10} \, \varepsilon^8 \Big ),\] 
where $I(n)$ is a positive constant (see \cite{Brendle4}, Proposition 19). Since $n \geq 52$, the function $\varepsilon \mapsto F(0,\varepsilon)$ has a strict local minimum at some point $\varepsilon_*$. The number $\varepsilon_*$ can be characterized as the unique positive solution of the equation 
\[\bigg ( 3 + \sqrt{9 - \frac{8(n+8)(n-8)}{(n+4)(n-10)}} \bigg ) \, \varepsilon_*^2 = \frac{2(n-8)}{n+4}.\] 
Moreover, by examining the Hessian of the function $F(\xi,\varepsilon)$ at the point $(0,\varepsilon_*)$,  one can show out that the function $(\xi,\varepsilon) \mapsto F(\xi,\varepsilon)$ attains a strict local minimum at the point $(0,\varepsilon_*)$. Hence, if $\rho^{2-n} \, \mu^{-2} \, \lambda^{n-10}$ is sufficiently small, then the function $(\xi,\varepsilon) \mapsto \lambda^{-8} \, \mu^{-2} \, \mathcal{F}_g(\lambda \xi,\lambda \varepsilon)$ has a local minimum at some nearby point $(\bar{\xi},\bar{\varepsilon})$. Consequently, the point $(\lambda \bar{\xi},\lambda \bar{\varepsilon})$ is a critical point of the function $(\xi,\varepsilon) \mapsto \mathcal{F}_g(\xi,\varepsilon)$. By Proposition \ref{lyapunov.schmidt}, the function $v := v_{(\lambda \bar{\xi},\lambda \bar{\varepsilon})}$ is a solution of the differential equation 
\[\Delta_g v - \frac{n-2}{4(n-1)} \, R_g \, v + n(n-2) \, |v|^{\frac{4}{n-2}} \, v = 0.\] 
Moreover, it can be shown that the function $v$ is nonnegative and 
\[\sup_{|x| \leq \lambda} v(x) \geq c(n) \, \lambda^{\frac{2-n}{2}}.\] From this, Theorem \ref{gluing} follows.

In \cite{Brendle-Marques}, it was shown that the Compactness Conjecture fails in dimension $25 \leq n \leq 51$.

\begin{theorem}[S.~Brendle, F.C.~Marques \cite{Brendle-Marques}]
\label{blow.up2}
Assume that $n \geq 25$. Then there exists a smooth Riemannian metric $g$ on $S^n$ and a sequence of positive functions $\{v_\nu: \nu \in \mathbb{N}\}$ with the following properties: 
\begin{itemize}
\item $g$ is not conformally flat. 
\item For each $\nu \in \mathbb{N}$, the function $v_\nu$ is a solution of the Yamabe equation 
\[\Delta_g v_\nu - \frac{n-2}{4(n-1)} \, R_g \, v_\nu + n(n-2) \, v_\nu^{\frac{n+2}{n-2}} = 0.\] 
\item For each $\nu \in \mathbb{N}$, we have $E_g(v_\nu) < Y(S^n)$. Moreover, $E_g(v_\nu) \to Y(S^n)$ as $\nu \to \infty$. 
\item $\sup_{S^n} v_\nu \to \infty$ as $\nu \to \infty$.
\end{itemize}
\end{theorem}

Note that the dimensional restriction in Theorem \ref{blow.up2} is sharp. The proof of Theorem \ref{blow.up2} is similar in spirit to the proof of Theorem \ref{blow.up}, though various adjustments and refinements are required. 

Using similar techniques, it is possible to construct counterexamples to the Weyl Vanishing Conjecture (as stated in the previous section); see \cite{Marques2} for details.

\section{A compactness result in dimension $n \leq 24$}

\label{compactness}

In this section, we describe the compactness result proved in \cite{Khuri-Marques-Schoen}. In order to state this result, we fix a compact manifold $(M,g_0)$. For any $p \in [1,\frac{n+2}{n-2}]$ we define
\begin{equation*}
\Phi_{p}=\{u > 0: \text{\rm $L_{g_0}u + n(n-2) \, u^p = 0$ on $M$}\},
\end{equation*}
where $L_{g_0} = \Delta_{g_0} -\frac{n-2}{4(n-1)} \, R_{g_0}$ denotes the conformal Laplacian. Although the geometric problem corresponds to the exponent $p = \frac{n+2}{n-2}$, which is critical with respect to the Sobolev embeddings, the consideration of the subcritical solutions is useful for the  purposes of developing a Morse theory and computing the total Leray-Schauder degree of the problem.

The main theorem of \cite{Khuri-Marques-Schoen} is:

\begin{theorem}[M.~Khuri, F.C.~Marques, R.~Schoen \cite{Khuri-Marques-Schoen}]
\label{compactness_thm}
Let $(M,g_0)$ be a compact Riemannian manifold of dimension $n\leq 24$, not conformally diffeomorphic to the standard sphere. Suppose $M$ is spin if $n\geq 8$. Then for any $\varepsilon \in (0,\frac{4}{n-2})$ there exists a constant $C>0$ depending only on $g_0$ and $\varepsilon$ such that
\begin{equation*}
C^{-1} \leq u \leq C \qquad \text{and} \qquad \|u\|_{C^2} \leq C,
\end{equation*}
for all $u \in \cup_{1+\varepsilon \leq p \leq \frac{n+2}{n-2}} \Phi_{p}$.
\end{theorem}

Let us now sketch the proof of Theorem \ref{compactness_thm}.

Similarly as before, we assume by contradiction that there exists a sequence $u_\nu \in \Phi_{p_\nu}$, with $p_\nu \in [1+\varepsilon,\frac{n+2}{n-2}]$, and such that $\sup_M u_\nu = u_\nu(y_\nu) \rightarrow \infty$ as $\nu \rightarrow \infty$. Standard blow-up analysis gives that $\lim_{\nu \rightarrow \infty} p_\nu = \frac{n+2}{n-2}$. Suppose that $\overline{y}$ is an isolated simple blow-up point. The first step is to obtain sharp approximations of the $u_\nu$ in a neighborhood of fixed size of $y_\nu$. These come from optimal pointwise estimates that generalize the estimates obtained by the second author in \cite{Marques1}. In high dimensions it is necessary to perform a refined blow-up analysis and go beyond the rotationally symmetric first approximation (standard bubble). The approximate solutions coincide with the ones introduced by the first author in \cite{Brendle3} to generalize the results of Aubin \cite{Aubin} and of Hebey and Vaugon \cite{Hebey-Vaugon}.
 
If $(x_1,\dots,x_n)$ are normal coordinates centered at $y_\nu$, we can view the metric as a smooth function taking values in the space of symmetric $n \times n$ matrices. It is convenient to work in conformal normal coordinates to simplify the computations. After performing a conformal change of the background metric $g_0$, one can assume that $\det g_0(x) = 1 + O(|x|^{2d+2})$ in geodesic normal coordinates around the point $y_\nu$ (cf. \cite{Lee-Parker}). We can write the background metric $g_0$ in the form 
\[g_0(x) = \exp(h(x))\] 
for some symmetric matrix $h_{ij}(x)$. As in \cite{Brendle3}, one looks at the Taylor expansion of $h_{ij}$ around the origin:
\[h_{ij}(x) = H_{ij}(x) + O(|x|^{d+1}),\] 
where $d=[\frac{n-2}{2}]$ and $H_{ij}(x)$ is a symmetric matrix whose entries are polynomials of degree less than or equal to $d$. Moreover, $H_{ij}(x)$ satisfies 
\begin{itemize}
\item $\sum_k H_{kk}(x)=0$,
\item $\sum_k x_k \, H_{ik}(x) = 0$,
\end{itemize}
for all $1 \leq i,j \leq n$ and $x \in \mathbb{R}^n$. Let us denote the vector space of such matrix-valued functions $H_{ij}(x)$ by $\mathcal{V}_n$. Note that $\mathcal{V}_n$ is finite-dimensional.

The optimal pointwise estimates established in \cite{Khuri-Marques-Schoen} can be used to expand the interior integral of (\ref{pohozaev_i}) in powers of $\varepsilon_\nu = u_\nu(y_\nu)^{-\frac{2}{n-2}}$. It turns out that the relevant correction terms are encoded in a canonical quadratic form $\mathcal{P}_n$ defined on $\mathcal{V}_n$. 

It is not difficult to check that the boundary integrals of (\ref{pohozaev_i}) are bounded by $C \, \varepsilon_\nu^{n-2}$, for some constant $C>0$. If $\mathcal{P}_n$ is a positive definite quadratic form on $\mathcal{V}_n$, then the Pohozaev identity (\ref{pohozaev_i}) implies that 
\[\sum_{s=2}^d \|H^{(s)}\|_{\mathcal{V}_n}^2 \, \varepsilon_\nu^{2s} \, |\log \varepsilon_\nu|^{\theta_s} \leq C \, \varepsilon_\nu^{n-2}.\] 
Here, $H^{(s)}$ denotes the homogeneous component of $H$ of degree $s$, $\|\cdot\|_{\mathcal{V}_n}$ is an arbitrary norm on the vector space $\mathcal{V}_n$, and $\theta_s$ is defined by 
\[\theta_s = \begin{cases} 1 & \text{\rm if $s = \frac{n-2}{2}$} \\ 0 & \text{\rm otherwise.} \end{cases}\] 
By taking the limit as $\nu \rightarrow \infty$, we conclude that $H_{ij}$ vanishes  at $\overline{y}$. This statement is equivalent to the Weyl Vanishing Conjecture. The sign contradiction with the Positive Mass Theorem follows similarly. 

Therefore, the problem is reduced to checking the positivity of $\mathcal{P}_n$. If $n$ is odd, for instance, and $\sum_{i,j}\partial_i \partial_j H_{ij} = 0$, the quadratic form is given by
\[\mathcal{P}_n(H,H) = \sum_{i,j,l} \sum_{s,t=2}^d c_{s+t} \int_{S_1^{n-1}(0)} \Big ( -\frac{1}{2} \partial_j H_{ij}^{(s)}\, \partial_l H_{il}^{(t)} + \frac{1}{4} \, \partial_l H_{ij}^{(s)} \, \partial_l H_{ij}^{(t)} \Big ).\] 
Here,
\[c_k = \int_0^\infty \frac{(r^2-1) \, r^{k+n-3}}{(1+r^2)^{n-1}} \, dr,\] 
for $k < n-2$. We refer the reader to the appendix of \cite{Khuri-Marques-Schoen} for a complete definition of $\mathcal{P}_n$.  

In \cite{Khuri-Marques-Schoen}, the eigenvalues of $\mathcal{P}_n$ are analyzed and the following result is proved: 

\begin{proposition}[\cite{Khuri-Marques-Schoen}]
The quadratic form $\mathcal{P}_n$, defined on $\mathcal{V}_n$, is positive definite if $n \leq 24$. Moreover, it has negative eigenvalues if $n \geq 25$.
\end{proposition}

We now describe some corollaries of Theorem \ref{compactness_thm}. In the case that all solutions of the Yamabe problem are nondegenerate, as will be the case for a generic conformal class
of Riemannian metrics, our previous results assert that there will be a finite number of solutions to the variational problem. Furthermore, denote by $N_\mu$ the number of solutions of (\ref{yamabe.equation.2}) of Morse index $\mu$, then we have he strong Morse inequalities  
\begin{equation*}
\sum_{\mu=0}^\lambda (-1)^{\lambda-\mu} \, N_\mu \geq (-1)^{\lambda}, \qquad \lambda=0,1,2,\ldots.
\end{equation*}
Indeed, it is well-known that the strong Morse inequalities hold for the sub-critical equation $L_{g_0} u + n(n-2) \, u^p = 0$ with $p < \frac{n+2}{n-2}$. Theorem \ref{compactness_thm} guarantees that the solutions of the sub-critical problem stay in a compact region as $p \to \frac{n+2}{n-2}$.

We also obtain:

\begin{corollary}[\cite{Khuri-Marques-Schoen}]
\label{intr_index}
Suppose that $(M,g)$ satisfies the assumptions of Theorem \ref{compactness_thm}, and assume that all critical points in $[g]$ are nondegenerate. Then there are a finite number of critical points $g_1,\ldots,g_k$ and we have
\begin{equation*}
\sum_{j=1}^k (-1)^{\text{\rm ind}(g_j)} = 1,
\end{equation*}
where $\text{\rm ind}(g_j)$ denotes the Morse index of the variational problem with volume constraint.
\end{corollary}

\section{The parabolic Yamabe flow}

In this final section, we discuss the Yamabe flow. More precisely, we fix a compact manifold $M$ and an initial metric $g_0$. We then evolve the metric by the evolution equation 
\begin{equation} 
\label{yamabe.flow}
\frac{\partial}{\partial t} g(t) = -(R_{g(t)} - r_{g(t)}) \, g(t), 
\end{equation}
where $r_{g(t)}$ denotes the mean value of the scalar curvature of $g(t)$: 
\[r_{g(t)} = \frac{\int_M R_{g(t)} \, d\text{\rm vol}_{g(t)}}{\text{\rm vol}(M,g(t))}.\] 
It is clear that the flow (\ref{yamabe.flow}) preserves the conformal structure. Hence, we may write $g(t) = u(t)^{\frac{4}{n-2}} \, g_0$, where $g_0$ denotes the initial metric. The conformal factor satisfies the evolution equation 
\begin{equation} 
\label{yamabe.flow2}
\frac{\partial}{\partial t} (u(t)^{\frac{n+2}{n-2}}) = \frac{n+2}{4} \, \Big ( \frac{4(n-1)}{n-2} \, \Delta_{g_0} u(t) - R_{g_0} \, u(t) + r_{g(t)} \, u(t)^{\frac{n+2}{n-2}} \Big ). 
\end{equation} 
The equation (\ref{yamabe.flow2}) is clearly parabolic, so the existence of a shorttime solution is guaranteed. In the late 1980s, Hamilton showed that, for any choice of the initial metric $g_0$, the Yamabe flow has a solution which exists for all $t \geq 0$. In other words, the Yamabe flow cannot develop a singularity in finite time. Inspired by this result, Hamilton proposed the following conjecture:

\begin{C2}[R.~Hamilton]
Let $g(t)$, $t \geq 0$, be a solution to the Yamabe flow on a compact manifold $M$. Then the metrics $g(t)$ converge to a metric of constant scalar curvature as $t \to \infty$.
\end{C2}

If the initial manifold $(M,g_0)$ is locally conformally flat and has positive Ricci curvature, then a result of Chow \cite{Chow} implies that the Yamabe flow converges to a metric of constant scalar curvature as $t \to \infty$. Using the method of moving planes, Ye \cite{Ye} proved the convergence of the Yamabe flow under the assumption that $(M,g_0)$ is locally conformally flat.

Schwetlick and Struwe \cite{Schwetlick-Struwe} proved the following convergence result: 

\begin{theorem}[H.~Schwetlick, M.~Struwe \cite{Schwetlick-Struwe}]
Let $(M,g_0)$ be a compact Riemannian manifold of dimension $n$, where $3 \leq n \leq 5$. Moreover, suppose that the Yamabe energy of $g_0$ is less than $\big [ Y(M,g_0)^{\frac{n}{2}} + Y(S^n)^{\frac{n}{2}} \big ]^{\frac{2}{n}}$. Finally, let $g(t)$, $t \geq 0$, be the unique solution to the Yamabe flow with initial metric $g_0$. Then the metrics $g(t)$ converge to a metric of constant scalar curvature as $t \to \infty$.
\end{theorem}

In 2005, the first author proved a convergence result in dimension $3 \leq n \leq 5$, which does not require any restrictions on the energy of the initial metric. This method also gives a new proof of Ye's convergence result in the locally conformally flat case.

\begin{theorem}[S.~Brendle \cite{Brendle1}]
\label{convergence}
Let $(M,g_0)$ be a compact Riemannian manifold of dimension $n$. We assume that either $3 \leq n \leq 5$ or $(M,g_0)$ is locally conformally flat. Moreover, let $g(t)$, $t \geq 0$, be the unique solution to the Yamabe flow with initial metric $g_0$. Then the metrics $g(t)$ converge to a metric of constant scalar curvature as $t \to \infty$.
\end{theorem}

In the remainder of this section, we will give an outline of the proof of Theorem \ref{convergence}. Note that the convergence of the flow follows directly from the maximum principle when $Y(M,g_0) \leq 0$. It therefore suffices to consider the case when $Y(M,g_0)$ is positive. The following result is the key ingredient in the proof of Theorem \ref{convergence}:

\begin{theorem}[S.~Brendle \cite{Brendle1}]
\label{key.ingredient}
Let $(M,g_0)$ be a compact Riemannian manifold of dimension $n$. We assume that either $3 \leq n \leq 5$ or $(M,g_0)$ is locally conformally flat. Moreover, suppose that $\kappa$ is a positive constant and $\{u_\nu: \nu \in \mathbb{N}\}$ is a sequence of positive functions satisfying 
\[\int_M u_\nu^{\frac{2n}{n-2}} \, d\text{\rm vol}_{g_0} = 1\] 
for all $\nu \in \mathbb{N}$ and 
\[\int_M \Big | \frac{4(n-1)}{n-2} \, \Delta_{g_0} u_\nu - R_{g_0} \, u_\nu + \kappa \, u_\nu^{\frac{n+2}{n-2}} \Big |^{\frac{2n}{n+2}} \, d\text{\rm vol}_{g_0} \to 0\] 
as $\nu \to \infty$. Then there exists a real number $\gamma > 0$ such that 
\begin{align*} 
&E_{g_0}(u_\nu) - \kappa \\ 
&\leq \bigg ( \int_M \Big | \frac{4(n-1)}{n-2} \, \Delta_{g_0} u_\nu - R_{g_0} \, u_\nu + \kappa \, u_\nu^{\frac{n+2}{n-2}} \Big |^{\frac{2n}{n+2}} \, d\text{\rm vol}_{g_0} \bigg )^{\frac{n+2}{2n} \, (1+\gamma)} 
\end{align*}
for $\nu$ sufficiently large.
\end{theorem}

The proof of Theorem \ref{key.ingredient} is quite involved. The argument can be simplified substantially if $(M,g_0)$ is conformally equivalent to the round sphere; this special case is discussed in more detail in \cite{Brendle2}. 

We now describe the main steps involved in the proof of Theorem \ref{key.ingredient}. After passing to a subsequence if necessary, the sequence $\{u_\nu: \nu \in \mathbb{N}\}$ converges to a function $u_\infty$ weakly in $H^1$. By a theorem of Struwe \cite{Struwe}, there exists an integer $m \geq 0$, a collection of positive real numbers $\{\varepsilon_{i,\nu}: 1 \leq i \leq m, \, \nu \in \mathbb{N}\}$, and a collection of points $\{y_{i,\nu}: 1 \leq i \leq m, \, \nu \in \mathbb{N}\} \subset M$ such that 
\[u_\nu - \sum_{i=1}^m \bigg ( \frac{4n(n-1)}{\kappa} \bigg )^{\frac{n-2}{4}} \, \Big ( \frac{\varepsilon_{i,\nu}}{\varepsilon_{i,\nu}^2 + d(y_{i,\nu},\cdot)^2} \Big )^{\frac{n-2}{2}} \to u_\infty\] 
strongly in $H^1$. There are two cases now: 

\textit{Case 1:} Suppose that the weak limit $u_\infty$ vanishes identically. In this case, we have 
\[\kappa = \lim_{\nu \to \infty} E_{g_0}(u_\nu) = \big [ m \, Y(S^n)^{\frac{n}{2}} \big ]^{\frac{2}{n}}.\] 
Moreover, one can show that 
\begin{align*} 
&E_{g_0}(u_\nu) - \kappa \\ 
&\leq C \, \bigg ( \int_M \Big | \frac{4(n-1)}{n-2} \, \Delta_{g_0} u_\nu - R_{g_0} \, u_\nu + \kappa \, u_\nu^{\frac{n+2}{n-2}} \Big |^{\frac{2n}{n+2}} \, d\text{\rm vol}_{g_0} \bigg )^{\frac{n+2}{n}} 
\end{align*}
for some uniform constant $C$ (see \cite{Brendle1}, Section 5).

\textit{Case 2:} Suppose that the weak limit $u_\infty$ does not vanish identically. In this case, 
\[\kappa = \lim_{\nu \to \infty} E_{g_0}(u_\nu) = \big [ E_{g_0}(u_\infty)^{\frac{n}{2}} + m \, Y(S^n)^{\frac{n}{2}} \big ]^{\frac{2}{n}}.\] 
Moreover, the function $u_\infty$ is a positive smooth solution of the partial differential equation 
\[\frac{4(n-1)}{n-2} \, \Delta_{g_0} u_\infty - R_{g_0} \, u_\infty + \kappa \, u_\infty^{\frac{n+2}{n-2}} = 0.\] 
In this case, we can show that 
\begin{align*} 
&E_{g_0}(u_\nu) - \kappa \\ 
&\leq C \, \bigg ( \int_M \Big | \frac{4(n-1)}{n-2} \, \Delta_{g_0} u_\nu - R_{g_0} \, u_\nu + \kappa \, u_\nu^{\frac{n+2}{n-2}} \Big |^{\frac{2n}{n+2}} \, d\text{\rm vol}_{g_0} \bigg )^{\frac{n+2}{2n} \, (1+\gamma)} 
\end{align*} 
for some positive constant $\gamma$ and some constant $C$. The proof uses a deep inequality for real analytic functions, which is due to Lojasiewicz (cf. \cite{Simon}). The details are presented in \cite{Brendle1}, Section 6. From this, Theorem \ref{key.ingredient} follows. \\

We now explain how Theorem \ref{convergence} follows from Theorem \ref{key.ingredient}. To that end, we consider a one-parameter family of metrics $g(t)$, $t \geq 0$, which evolves by the Yamabe flow. 

\begin{proposition}[\cite{Brendle1}]
\label{integrability}
Let $(M,g_0)$ be a compact Riemannian manifold. We assume that either $3 \leq n \leq 5$ or $(M,g_0)$ is locally conformally flat. Moreover, let $g(t)$, $t \geq 0$, be the unique solution to the Yamabe flow with initial metric $g_0$.
Then 
\[\int_0^\infty \bigg ( \int_M (R_{g(\tau)} - r_{g(\tau)})^2 \, d\text{\rm vol}_{g(\tau)} \bigg )^{\frac{1}{2}} \, d\tau < \infty.\] 
\end{proposition} 

\textit{Sketch of the proof of Proposition \ref{integrability}.} Without loss of generality, we may assume that $\text{\rm vol}(M,g(t)) = 1$ for all $t \geq 0$. As above, we denote by $r_{g(t)}$ the mean value of the scalar curvature of $g(t)$. Then the function $t \mapsto r_{g(t)}$ is decreasing; in particular, the limit $\kappa := \lim_{t \to \infty} r_{g(t)}$ exists. Moreover, one can show that 
\[\int_M |R_{g(t)} - \kappa|^{\frac{2n}{n+2}} \, d\text{\rm vol}_{g(t)} \to 0\] 
as $t \to \infty$. We now write $g(t) = u(t)^{\frac{4}{n-2}} \, g_0$ for some positive function $u(t)$. Then 
\[\int_M u(t)^{\frac{2n}{n-2}} \, d\text{\rm vol}_{g_0} = 1\] 
for all $t \geq 0$ and 
\[\int_M \Big | \frac{4(n-1)}{n-2} \, \Delta_{g_0} u(t) - R_{g_0} \, u(t) + \kappa \, u(t)^{\frac{n+2}{n-2}} \Big |^{\frac{2n}{n+2}} \, d\text{\rm vol}_{g_0} \to 0\] 
as $t \to \infty$. We claim that there exists a positive constant $\gamma$ such that 
\begin{align} 
\label{estimate}
&E_{g_0}(u(t)) - \kappa \notag \\ 
&\leq \bigg ( \int_M \Big | \frac{4(n-1)}{n-2} \, \Delta_{g_0} u(t) - R_{g_0} \, u(t) + \kappa \, u(t)^{\frac{n+2}{n-2}} \Big |^{\frac{2n}{n+2}} \, d\text{\rm vol}_{g_0} \bigg )^{\frac{n+2}{2n} \, (1+\gamma)} 
\end{align}
for $t$ sufficiently large. Indeed, if (\ref{estimate}) is false, we can find a sequence of times $t_\nu$ such that $t_\nu \geq \nu$ and 
\begin{align*} 
&E_{g_0}(u(t_\nu)) - \kappa \\ 
&> \bigg ( \int_M \Big | \frac{4(n-1)}{n-2} \, \Delta_{g_0} u(t_\nu) - R_{g_0} \, u(t_\nu) + \kappa \, u(t_\nu)^{\frac{n+2}{n-2}} \Big |^{\frac{2n}{n+2}} \, d\text{\rm vol}_{g_0} \bigg )^{\frac{n+2}{2n} \, (1+\frac{1}{\nu})} 
\end{align*}
for each $\nu$, and this contradicts Theorem \ref{key.ingredient}. Since $(M,g(t))$ has unit volume, the estimate (\ref{estimate}) implies 
\begin{align*} 
r_{g(t)} - \kappa 
&\leq \bigg ( \int_M |R_{g(t)} - \kappa|^{\frac{2n}{n+2}} \, d\text{\rm vol}_{g(t)} \bigg )^{\frac{n+2}{2n} \, (1+\gamma)} \\ 
&\leq \bigg ( \int_M (R_{g(t)} - \kappa)^2 \, d\text{\rm vol}_{g(t)} \bigg )^{\frac{1+\gamma}{2}} \\ 
&= \bigg ( \int_M (R_{g(t)} - r_{g(t)})^2 \, d\text{\rm vol}_{g(t)} + (r_{g(t)} - \kappa)^2 \bigg )^{\frac{1+\gamma}{2}}. 
\end{align*} 
Hence, if $t$ is sufficiently large, then we obtain 
\[r_{g(t)} - \kappa \leq C \, \bigg ( \int_M (R_{g(t)} - r_{g(t)})^2 \, d\text{\rm vol}_{g(t)} \bigg )^{\frac{1+\gamma}{2}}\] 
for some uniform constant $C$. On the other hand, we have 
\[r_{g(t)} - \kappa = \frac{n-2}{2} \int_t^\infty \int_M (R_{g(\tau)} - r_{g(\tau)})^2 \, d\text{\rm vol}_{g(\tau)} \, d\tau.\] 
A standard ODE lemma now implies that 
\[\int_0^\infty \bigg ( \int_M (R_{g(\tau)} - r_{g(\tau)})^2 \, d\text{\rm vol}_{g(\tau)} \bigg )^{\frac{1}{2}} \, d\tau < \infty,\] 
as claimed. \\

Using Proposition \ref{integrability}, we are able to rule out volume concentration: 

\begin{proposition}[\cite{Brendle1}]
\label{volume.concentration}
Let $(M,g_0)$ be a compact Riemannian manifold. We assume that either $3 \leq n \leq 5$ or $(M,g_0)$ is locally conformally flat. Moreover, let $g(t)$, $t \geq 0$, be the unique solution to the Yamabe flow with initial metric $g_0$. Then, given any positive real number $\eta$, we can find a real number $r > 0$ such that 
\[\text{\rm vol}(B_r(p),g(t)) \leq \eta\] 
for all points $p \in M$ and all $t \geq 0$. Here, $B_r(p)$ denotes a geodesic ball of radius $r$ with respect to the background metric $g_0$. 
\end{proposition}

\textit{Sketch of the proof of Proposition \ref{volume.concentration}.} It follows from Proposition \ref{integrability} that 
\[\int_0^\infty \int_M |R_{g(\tau)} - r_{g(\tau)}| \, d\text{\rm vol}_{g(\tau)} \, d\tau < \infty.\] 
Consequently, we can find a real number $T > 0$ such that 
\[\int_T^\infty \int_M |R_{g(\tau)} - r_{g(\tau)}| \, d\text{\rm vol}_{g(\tau)} \, d\tau \leq \frac{\eta}{n}.\] 
We next choose a real number $r > 0$ such that $\text{\rm vol}(B_r(p),g(t)) \leq \frac{\eta}{2}$ for all points $p \in M$ and all $t \in [0,T]$. Then 
\begin{align*} 
\text{\rm vol}(B_r(p),g(t)) 
&= \text{\rm vol}(B_r(p),g(T)) - \frac{n}{2} \int_T^t \int_{B_r(p)} (R_{g(\tau)} - r_{g(\tau)}) \, d\text{\rm vol}_{g(\tau)} \, d\tau \\ 
&\leq \text{\rm vol}(B_r(p),g(T)) + \frac{n}{2} \int_T^t \int_M |R_{g(\tau)} - r_{g(\tau)}| \, d\text{\rm vol}_{g(\tau)} \, d\tau \\ 
&\leq \eta 
\end{align*} 
for all points $p \in M$ and all $t \in [T,\infty)$. From this, Proposition \ref{volume.concentration} follows. \\

Once we know that volume concentration does not occur, it follows from standard arguments that the function $u(t)$ is uniformly bounded from above and below. From this, the convergence of the flow follows. This completes our sketch of the proof of Theorem \ref{convergence}.

Finally, we state a generalization of Theorem \ref{convergence} to the higher dimensional setting. Let $(M,g_0)$ be a compact Riemannian manifold of dimension $n \geq 6$, and let $d = [\frac{n-2}{2}]$. We denote by $\mathcal{Z}$ the set of all points $\overline{y} \in M$ such that 
\[\limsup_{y \to \overline{y}} d(\overline{y},y)^{2-d} \, |W_{g_0}(y)| = 0,\] 
where $W_{g_0}$ denotes the Weyl tensor of the background metric $g_0$. Note that the set $\mathcal{Z}$ depends only on the conformal class of $g_0$.

\begin{theorem}[S.~Brendle \cite{Brendle3}]
\label{higher.dim}
Let $(M,g_0)$ be a compact Riemannian manifold of dimension $n \geq 6$. We assume that either 
$M$ is spin or $\mathcal{Z} = \emptyset$. Moreover, let $g(t)$, $t \geq 0$, be the unique solution to the Yamabe flow with initial metric $g_0$. Then the metrics $g(t)$ converge to a metric of constant scalar curvature as $t \to \infty$.
\end{theorem}

The proof of Theorem \ref{higher.dim} is similar in spirit to that of Theorem \ref{convergence}. In order to extend Proposition \ref{key.ingredient} to higher dimensions, one needs to construct a suitable family of test functions with Yamabe energy less than $Y(S^n)$. These test functions are constructed by a generalization of Aubin's method (cf. \cite{Aubin}, \cite{Hebey-Vaugon}); see \cite{Brendle3} for details.

\end{document}